\documentclass[11pt]{elsart}

\usepackage{url}

\usepackage{amssymb}
\usepackage{times}
\usepackage{enumerate}
{\bfseries}{\itshape}
{\bfseries}{\itshape}

\newcommand {\bbox}{\rule{0.6em}{0.6em}}

\newcommand{\boxi}{\ensuremath{{box}}}

\newcommand{\tw}{\ensuremath{{tw}}}

\newcommand{\ignore}[1]{}

\renewenvironment{pf}
{\begin{oldproof}}{\hfill\bbox\end{oldproof}}

\begin{document}

\begin{frontmatter}

\title{Boxicity and Maximum degree}

\author[sun]{L. Sunil Chandran} 
\and
\author[sun]{Mathew C. Francis} 
\and 
\author[ns]{Naveen Sivadasan} 
\address[sun]{Indian Institute of Science,
Dept. of Computer Science and Automation,
Bangalore  560012, India.  email: \emph{sunil,mathew@csa.iisc.ernet.in}}
\address[ns]{Strand Life Sciences,
237, Sir. C. V. Raman Avenue, Rajmahal Vilas,
Bangalore  560080, India.  email: \emph{naveen@strandls.com}}

\date{}
\bibliographystyle{plain}

\pagenumbering{arabic}
\maketitle

\begin {abstract}

An axis-parallel  $d$--dimensional box is a Cartesian product 
$R_1 \times R_2 \times \cdots \times R_d$ where $R_i$ (for $1 \le i \le d$)
is a closed interval of the form $[a_i, b_i]$ on the real line.
For a  graph $G$, its \emph{boxicity} $\boxi(G)$ is the minimum dimension $d$, such that
$G$ is representable as the intersection graph of 
(axis--parallel) boxes in $d$--dimensional space. The concept of
boxicity finds applications in various areas such as ecology, operation research etc.

We show that for any graph $G$ with maximum degree $\Delta$, $\boxi(G) \le 2 \Delta^2 + 2$.
That the bound does not depend on the number of vertices
is a bit surprising considering the fact that there are highly connected bounded degree graphs such
as expander graphs. Our proof is very short and constructive. We conjecture that $\boxi(G)$ is $O(\Delta)$.

\end {abstract}

\end{frontmatter}

Let ${\mathcal{F}} = \{ S_x \subseteq U : x \in V \}$ be a family of 
subsets of a universe $U$, where $V$ is an index set. The intersection graph 
$\Omega({\mathcal{F}})$ of ${\mathcal{F}}$ 
has $V$ as vertex set, and two distinct vertices
$x$ and $y$ are adjacent if and only if $S_x \cap S_y \ne \emptyset$.

Representation of graphs as the intersection graphs of various
geometrical objects is a well studied topic in graph theory.  
A prime example of a graph class defined in this way is the class of 
interval graphs: {\it A graph $G$ is an interval graph \emph{if and only if} 
$G$ has an interval realization: i.e., each
vertex of  $G$ can be associated to an interval on the real 
line such that two intervals intersect if and only if
the corresponding vertices are adjacent.}
Motivated by theoretical as well as practical considerations, 
graph theorists have tried to generalize the concept of interval graphs
in various ways. 
One such generalization is the concept of \emph{boxicity} defined as follows.

An axis-parallel  $d$--dimensional box is a Cartesian product 
$R_1 \times R_2 \times \cdots \times R_d$ where $R_i$ (for $1 \le i \le d$)
is a closed interval of the form $[a_i, b_i]$ on the real line.
For a  graph $G$, its \emph{boxicity} $\boxi(G)$ is the minimum dimension $d$, such that
$G$ is representable as the intersection graph of 
(axis--parallel) boxes in $d$--dimensional space. It is easy
to see that the class of graphs with $d\le 1$ is exactly the class of
interval graphs.
\ignore{
	A $b$--dimensional \emph{box representation} of a graph $G=(V,E)$ is
	a mapping  that maps each
	each $u \in V$ to an axis-parallel $b$--dimensional box $B_u$
	such that $G$ is the intersection graph of the family $\{B_u : u \in V\}$.
}

The concept of boxicity was  introduced by F. S.  Roberts \cite
{Roberts}. 
It  finds applications in niche overlap (competition) in ecology and to
problems of fleet maintenance in operations research. (See \cite {CozRob}.) 
It was shown 
by Cozzens \cite {Coz}  that computing the boxicity of a graph is NP--hard. 
This
was later strengthened by Yannakakis \cite {Yan1}, 
and finally by  Kratochvil \cite {Kratochvil} 
who  showed that deciding  whether
boxicity of a graph is at most 2 itself 
is NP--complete. 

There have  been many attempts to estimate or bound the boxicity of 
graph classes  with special structure. In his pioneering work,
F. S. Roberts proved that the boxicity of complete $k$--partite graphs
are $k$.  Scheinerman \cite {Scheiner} 
showed that the boxicity of outer planar graphs is at most $2$.
Thomassen \cite {Thoma1} proved that 
the boxicity of planar graphs is
bounded above by $3$. The boxicity of split graphs is investigated by
Cozzens and Roberts \cite{CozRob}. In a recent manuscript \cite{CN05}
the authors showed that $\boxi(G) \le \tw(G) + 2$ where
$\tw(G)$ is the treewidth of $G$.
Little is known about the structure imposed on a graph by its high boxicity.

A number of NP-hard problems are known to be polynomial time solvable
for interval graphs. Since boxicity is a direct generalization of the notion of 
interval graphs, such results may generalize to bounded boxicity graphs.
Thus our result may be of interest from an algorithmic point of view.

Researchers have also tried to generalize or extend  the
concept of boxicity in various ways. The poset boxicity \cite {TroWest}, 
the rectangular number \cite {ChangWest}, grid dimension \cite {Bellantoni},
circular dimension \cite {Feinberg,Shearer}  and the boxicity 
of digraphs \cite {ChangWest1} are some  examples. 

Let $G$ be a simple, finite, undirected, unweighted graph on $n$ vertices. 
Let $V(G)$ denote the vertex set of $G$ and $E(G)$ denote the edge set of $G$.
Let $\Delta(G)$ denote the maximum degree of $G$.
Let $I_1, \ldots, I_k$ be $k$
interval graphs such that $V(I_j) = V(G)$ for $1 \le j \le k$.
If $E(G) = E(I_1) \cap \cdots \cap E(I_k),$ then we say that
$I_1, \ldots, I_k$ is an \emph{interval graph representation} of $G$.
The following equivalence is well-known.

\noindent
\textbf{Fact [Roberts \cite{Roberts}]}\emph{~The minimum $k$ such that there exists an interval graph representation
of $G$ using $k$ interval graphs $I_1, \ldots, I_k$  is the same
as $\boxi(G)$.}


\begin{thm} 
For any graph $G$ with maximum degree $\Delta$, $\boxi(G) \le 2\Delta^2 + 2$. 
\end{thm}

\begin{pf}
Let $V(G) = V$.
Let $G^2$ denote the \emph{square} of $G$ defined as follows: $V(G^2) = V$
and two vertices $u$ and $v$ are adjacent in $G^2$ if and only if
the shortest distance between $u$ and $v$ in $G$ is either $1$ or $2$.
Let $\chi(G^2) = k$, where $\chi(H)$ denote the \emph{chromatic number} of the graph $H$.

Consider an optimal vertex coloring $c: V \rightarrow \{1,\ldots,k\}$ of $G^2$.
For $1 \le i \le k$, let $V_i = \{ u \in V ~|~ c(u) = i \}$ be the $i$th
color class. 
For $1 \le i \le k$, let $G_i$ be defined as follows. $V(G_i) = V$ and
$E(G_i) = E(G) \cup \{ (u, v): u, v \in V - V_i \mbox{~and~} u \not= v\}$. 
We claim that  $E(G) = E(G_1) \cap \cdots \cap E(G_k)$. To see this,
first observe that for $1 \le i \le k$, $E(G) \subseteq E(G_i)$. Now consider $(u, v) \notin E(G)$.
Let $c(u) = i$. Then by construction, $(u, v) \notin E(G_i)$.

We now show that $\boxi(G_i) \le 2$ for $1 \le i \le k$.
First we claim that in $G$, for any vertex $w \in V - V_i$, 
$w$ has at most one neighbor in $V_i$. 
This is because, if $w$ has two neighbors say $x$ and $y$ in $V_i$, 
then clearly $x$ is adjacent
to $y$ in $G^2$ and thus they can not belong to the same color class $V_i$.
Now, by construction of $G_i$, we have for any $(u, v) \in V_i \times (V - V_i)$,
$(u, v)$ belongs to $E(G_i)$ if and only if $(u, v)$ belongs to $E(G)$. 
Thus it follows that with respect to $G_i$ also, for any vertex $w \in V - V_i$,
$w$ has at most one neighbor in $V_i$.
Without loss of generality, let the vertices in $V_i$ be $\{1, \ldots, h\}$ where
$h = |V_i|$. Consider two orderings $\pi_0$ and $\pi_1$ of $V_i$ such that
for any $j \in V_i$, $\pi_0(j) = j$ and $\pi_1(j) = h - j + 1$.
For $r \in \{0, 1\}$, define the interval graph $I_r$ on the vertex set $V$ as follows:
For $w \in V_i$, let the interval $[\pi_r(w), \pi_r(w)]$ be assigned to $w$.
For $w \in V - V_i$,
if $w$ has no neighbors in $V_i$ with respect to $G_i$ then assign the interval $[0, 0]$ to $w$.
Otherwise let $z$ be its only neighbor in $V_i$ with respect to $G_i$. 
 Assign the interval $[0, \pi_r(z)]$ to $w$. 
We claim that $E(G_i) = E(I_0) \cap E(I_1)$ and thus $\boxi(G_i) \le 2$.
By construction, it is clear that $E(G_i) \subseteq E(I_r)$ for $r \in \{0, 1\}$.
It remains to show that  if $(u, v) \notin E(G_i)$ then
either $(u, v) \notin E(I_0)$ or $(u, v) \notin E(I_1)$. Since $V - V_i$ induces
a complete graph in $G_i$, if $(u, v) \notin E(G_i)$ then either $u,v \in V_i$ or $u \in V_i$ and $v \in V - V_i$.
Clearly, by the construction, $V_i$ is an independent set in $I_0$ as well as $I_1$. Thus
the only case we have to consider is when $u \in V_i$ and $v \in V - V_i$.
If $v$ has no neighbors in $V_i$ then clearly $(u, v) \notin E(I_0)$ and $(u, v) \notin E(I_1)$.
Otherwise let $x$ be the (only) neighbor of $v$ in $V_i$. Now, clearly either $\pi_0(u) > \pi_0(x)$
or $\pi_1(u) > \pi_1(x)$. It follows that $(u, v) \notin E(I_0)$ or $(u, v) \notin E(I_1)$.

Recalling that $E(G) = \bigcap_{i=1}^k E(G_i)$,
it follows that  $\boxi(G) \le \sum_{i=1}^k \boxi(G_i) \le 2k$.
Now using the well-known fact that 
for any graph $H$, $\chi(H) \le \Delta(H) + 1$, (see chapter 5, \cite{Diest}) it follows that
$k \le \Delta(G^2) + 1 \le \Delta^2 + 1$ and the result follows.
\end{pf}

\noindent
\textbf{Remark.} We conjecture that $\boxi(G)$ is $O(\Delta)$. In fact, given any $\Delta$,
it is not difficult to construct graphs of boxicity $\Omega(\Delta)$ on arbitrarily large number of vertices,
using a construction given by Roberts \cite{Roberts}.



\begin{thebibliography}{10}

\bibitem{Bellantoni}
S.~Bellantoni, I.~Ben-Arroyo Hartman, T.~Przytycka, and S.~Whitesides.
\newblock Grid intersection graphs and boxicity.
\newblock {\em Discrete mathematics}, 114(1-3):41--49, April 1993.

\bibitem{CN05}
L.~Sunil Chandran and Naveen Sivadasan.
\newblock Treewidth and boxicity.
\newblock Submitted, Available at http://arxiv.org/abs/math.CO/0505544, 2005.

\bibitem{ChangWest1}
Y.~W. Chang and Douglas~B. West.
\newblock Interval number and boxicity of digraphs.
\newblock In {\em Proceedings of the 8th International Graph Theory Conf.},
  1998.

\bibitem{ChangWest}
Y.~W. Chang and Douglas~B. West.
\newblock Rectangle number for hyper cubes and complete multipartite graphs.
\newblock In {\em 29th {SE} conf. Comb., Graph Th. and Comp., Congr. Numer.
  132}, pages 19--28, 1998.

\bibitem{Coz}
M.~B. Cozzens.
\newblock Higher and multidimensional analogues of interval graphs.
\newblock Ph. D thesis, Rutgers University, New Brunswick, NJ, 1981.

\bibitem{CozRob}
M.~B. Cozzens and F.~S. Roberts.
\newblock Computing the boxicity of a graph by covering its complement by
  cointerval graphs.
\newblock {\em Discrete Applied Mathematics}, 6:217--228, 1983.

\bibitem{Diest}
Reinhard Diestel.
\newblock {\em Graph Theory}, volume 173.
\newblock Springer Verlag, New York, 2 edition, 2000.

\bibitem{Feinberg}
Robert~B. Feinberg.
\newblock The circular dimension of a graph.
\newblock {\em Discrete mathematics}, 25(1):27--31, 1979.

\bibitem{Kratochvil}
J.~Kratochvil.
\newblock A special planar satisfiability problem and a consequence of its
  {NP}--completeness.
\newblock {\em Discrete Applied Mathematics}, 52:233--252, 1994.

\bibitem{Roberts}
F.~S. Roberts.
\newblock {\em Recent Progresses in Combinatorics}, chapter On the boxicity and
  Cubicity of a graph, pages 301--310.
\newblock Academic Press, New York, 1969.

\bibitem{Scheiner}
E.~R. Scheinerman.
\newblock Intersectin classes and multiple intersection parameters.
\newblock Ph. D thesis, Princeton University, 1984.

\bibitem{Shearer}
J.~B. Shearer.
\newblock A note on circular dimension.
\newblock {\em Discrete mathematics}, 29(1):103--103, 1980.

\bibitem{Thoma1}
C.~Thomassen.
\newblock Interval representations of planar graphs.
\newblock {\em Journal of combinatorial theory, Ser {B}}, 40:9--20, 1986.

\bibitem{TroWest}
W.~T. Trotter and Jr. Douglas~B. West.
\newblock Poset boxicity of graphs.
\newblock {\em Discrete Mathematics}, 64(1):105--107, March 1987.

\bibitem{Yan1}
Mihalis Yannakakis.
\newblock The complexity of the partial order dimension problem.
\newblock {\em {SIAM} Journal on Algebraic Discrete Methods}, 3:351--358, 1982.

\end{thebibliography}

\end {document}